\documentclass[11pt,reqno]{amsart} 
\usepackage{amsmath}
\usepackage{amsthm}
\usepackage{amssymb}
\usepackage{amsfonts}
\usepackage{latexsym}

\usepackage{latexsym}
\usepackage{graphicx,color}
\usepackage{enumerate,fancybox}

\newtheorem{theorem}{Theorem}[section]

\begin{document}
\title[Fine asymptotics for Bergman polynomials]
{Fine asymptotics for Bergman polynomials over domains with corners}

\date{\today}

\author[N. Stylianopoulos]{Nikos Stylianopoulos}
\address{Department of Mathematics and Statistics,
University of Cyprus,
P.O. Box 20537,
1678 Nicosia,
Cyprus}
\email{nikos@ucy.ac.cy
}
\urladdr{http://ucy.ac.cy/\textasciitilde nikos}
\subjclass[2000]{}

\maketitle

\begin{abstract}
Let $G$ be a bounded simply-connected domain in the complex plane $\mathbb{C}$,
whose boundary $\Gamma:=\partial G$ is a Jordan curve, and let
$\{p_n\}_{n=0}^{\infty}$ denote the sequence of Bergman polynomials of
$G$. This is defined as the sequence
$$
p_n(z) = \lambda_n z^n+\cdots, \quad \lambda_n>0,\quad n=0,1,2,\ldots,
$$
of polynomials that are orthonormal with respect to the inner product
$$
\langle f,g\rangle := \int_G f(z) \overline{g(z)} dA(z),
$$
where $dA$ stands for the area measure.

The purpose of this note is to report on recent results regarding the fine asymptotic
behaviour of the the leading coefficients $\lambda_n$, $n\in\mathbb{N}$, and the polynomials $p_n(z)$, for $z\in\Omega:=\mathbb{C}\setminus(\overline{G})$,  in cases when the boundary $\Gamma$ contains corners. These results complement an investigation started in the 1920's by T.\ Carleman, who obtained the fine asymptotics for domains with analytic boundaries and carried over by P.K.\ Suetin in the 1960's, who established them for domains with smooth boundaries.
\end{abstract}

\section{Introduction}\label{intro}
Let $G$ be a bounded simply-connected domain in the complex
plane, and let $\{p_n\}_{n=0}^{\infty}$ denote the sequence of {\it
Bergman orthogonal polynomials} of $G$. These are defined as
the sequence
$$
p_n(z) = \lambda_n z^n+ \lambda_{n-1} z^{n-1} + \cdots, \quad
\lambda_n>0,\quad n=0,1,2,\ldots,
$$
of polynomials which are orthonormal with respect to the inner product
$$
\langle f,g\rangle := \int_G f(z) \overline{g(z)} dA(z),
$$
where $dA$ stands for the area  measure.

Let $\Omega:=\overline{\mathbb{C}}\setminus\overline{G}$ denote the complement (in $\overline{\mathbb{C}}$) of $\overline{G}$ and let $\Phi$ denote the conformal map  $\Omega\to\Delta:=\{w:|w|>1\}$, normalized so that, near infinity,
\begin{equation}\label{eq:Phi}
\Phi(z)=\gamma z+\gamma_0+\frac{\gamma_1}{z}+\frac{\gamma_2}{z^2}+\cdots,\quad \gamma>0.
\end{equation}
We note that $1/\gamma$ gives the (logarithmic) capacity $\textup{cap}(\Gamma)$ of $\Gamma:=\partial G$.

Our main task is to consider the asymptotic behaviour of the leading coefficients $\lambda_n$ and the Bergman polynomials $p_n(z)$, for $z\in\Omega$.

\section{Fine asymptotics and applications}
The first result regarding the fine asymptotics of $\lambda_n$ and $p_n$ was
derived by T.\ Carleman, for domains bounded by analytic Jordan curves. In this case the conformal map $\Phi$ has an analytic and one-to-one continuation across $\Gamma$ into $G$. Below we use $\Gamma_R$ to denote the level line  $\{z:|\Phi(z)|=R\}$ of the mapping $\Phi$.

\begin{theorem}\label{tm:Carleman} \textup{(\cite{Ca23})}
If $\rho<1 $ is the smallest index for which $\Phi$ is
conformal in $\textup{ext}(\Gamma_\rho)$, then
$$
{\frac{n+1}{\pi}\frac{\gamma^{2(n+1)}}{\lambda_n^2}=1-\alpha_n,}\,\,
\mbox{ where } 0\le \alpha_n\le c_1(\Gamma)\,
\rho^{2n}
$$
and
$$
{p_n(z)=\sqrt{\frac{n+1}{\pi}}\Phi^n(z)\Phi^\prime(z)\{1+A_n(z)\},}
\quad n\in\mathbb{N},\quad
$$
where
$$
|A_n(z)|\le c_2(\Gamma)\sqrt{n}\,\rho^n,\quad z\in\overline{\Omega}.
$$
\end{theorem}
(Here and in the sequel we use $c(\Gamma)$, $c_1(\Gamma)$ and $c_2(\Gamma)$ to denote  positive constants that depend only on $\Gamma$.)

The next available result is due to P.K.\ Suetin and require that the boundary curve $\Gamma$ belongs to a smoothness class. We say that $\Gamma\in C(p,\alpha)$, for some $p\in \mathbb{N}$ and $0<\alpha<1$, if $\Gamma$ is given by $z=g(s)$, where $s$ is the arclength, with $g^{(p)}\in \textup{Lip}\alpha$. Then, both $\Phi$ and $\Psi:=\Phi^{-1}$ are p times continuously differentiable in $\overline{\Omega}\setminus\{\infty\}$ and $\overline{\Delta}\setminus\{\infty\}$ respectively, with $\Phi^{(p)}$ and $\Psi^{(p)}\in\textup{Lip}\alpha$.
\begin{theorem}\label{th:Suetin}\textup{(\cite{Su74})}
Assume that $\Gamma\in C(p+1,\alpha)$, with $p+\alpha>1/2$. Then
$$
{\frac{n+1}{\pi}\frac{\gamma^{2(n+1)}}{\lambda_n^2}=1-\alpha_n,}\,\,\mbox{ where }
0\le \alpha_n\le c_1(\Gamma)\,
\frac{1}{n^{2(p+\alpha)}}
$$
and
$$
{p_n(z)=\sqrt{\frac{n+1}{\pi}}\Phi^n(z)\Phi^\prime(z)\{1+A_n(z)\},}\quad n\in\mathbb{N},\quad
$$
where
$$
|A_n(z)|\le c_2(\Gamma)\,\frac{\log n}{n^{p+\alpha}},\quad z\in\overline{\Omega}.
$$
\end{theorem}

Apart from the above two results, the only other, of this kind, which are known are due to E.R.\ Johnston. However, Johnston's results are derived under analytic assumptions on certain functions associated to the conformal maps $\Phi$ and $\Psi$ (as compared to the geometric assumptions on $\Gamma$ in Theorems~\ref{tm:Carleman} and \ref{th:Suetin}) and they do not provide the order of decay of the error terms $\alpha_n$ and $A_n(z)$. (For an account of these results we refer the reader to \cite{RW}.)

The lack of fine asymptotics for non-smooth boundaries, even for special cases, makes one to incline to think that they do no actually exists. This would have been the case with us, had we not have available plenty of numerical data. (This is essential, because we are not aware of a single case of non-smooth $\Gamma$ for which the leading coefficients $\lambda_n$, $n=0,1,\ldots$, are known explicitly in terms of $n$.) Table~\ref{tab:1} contains a range of computed values of the error term $\alpha_n$,  for the case of a semi-disk. They were obtained after constructing the Bergman polynomials by means of the Gram-Schmidt process, and using the known exact value of $\gamma$. Having computed the values of $\alpha_n$ (and in view of the results of Suetin) we test the hypothesis
$
\displaystyle{\alpha_n \approx\,C\,\frac{1}{n^s}.}
$

\begin{table}[ht]
\begin{center}
\renewcommand{\baselinestretch}{0.9}
\begin{tabular}{|c|cc|}
\hline $ $
$n$& $\alpha_n$           & $s$   \\ \hline
51 & 0.003\,263\,458\,678 &  -       \\
52 & 0.003\,200\,769\,764 & 0.998\,887     \\
53 & 0.003\,140\,444\,435 & 0.998\,899     \\
54 & 0.003\,082\,351\,464 & 0.998\,911     \\
55 & 0.003\,026\,369\,160 & 0.998\,923     \\
56 & 0.002\,972\,384\,524 & 0.998\,934     \\
57 & 0.002\,920\,292\,482 & 0.998\,946    \\
58 & 0.002\,869\,952\,027 & 0.998\,957    \\
59 & 0.002\,821\,401\,485 & 0.998\,968     \\
60 & 0.002\,774\,426\,207 & 0.998\,979   \\ \hline
\end{tabular}
\normalsize
\end{center}
\caption{The rate of decay of $\alpha_n$ for the unit semi-disk.} \label{tab:1}
\end{table}
The numbers in Table~\ref{tab:1} indicate clearly that
$
{\alpha_n \approx\, C\,\frac{1}{n}}.
$
The same behaviour was observed in a number of other different non-smooth cases, with various angles involved. Accordingly, we have made conjectures regarding fine asymptotics in \cite{Baetal} and \cite{NS04}.
The main the purpose of this note is to report on results that verify these conjectures.

To this end, we have the following:
\begin{theorem}\label{th:NS1}
Assume that $\Gamma$ is piecewise analytic without cusps. Then, 
$$
{\frac{n+1}{\pi}\frac{\gamma^{2(n+1)}}{\lambda_n^2}=1-\alpha_n},
$$
where
$$
0\le\alpha_n\le c(\Gamma)\,\frac{1}{n},\quad n\in\mathbb{N}.
$$
\end{theorem}

\begin{theorem}\label{th:NS2}
Assume that $\Gamma$ is piecewise analytic without cusps. Then, 
for any $z\in\Omega$,
$$
{p_n(z)=\sqrt{\frac{n+1}{\pi}}\Phi^n(z)\Phi^\prime(z)
\left\{1+A_n(z)\right\}},
$$
where
$$
|A_n(z)|\le \frac{c(\Gamma)}{\textup{dist}(z,\Gamma)|\Phi^\prime(z)|}\,\frac{1}{\sqrt{n}},\quad n\in\mathbb{N}.
$$ 
\end{theorem}

As applications, we report the following two theorems that refine well-known classical results, regarding the position of zeros and the weak (or $n$-th root) asymptotics of the sequence $\{p_n(z)\}_{n=1}^\infty$. The first complements Fejer's theorem which states that \textit{all the zeros of} $p_n(z)$, $n\in\mathbb{N}$, \textit{are contained on the convex hull} of $\overline{G}$. The second extends the applicability of an important theorem of Stall and Totik, cf. \cite[Thm 3.3.1(ii)]{StTobo}.
\begin{theorem}
Assume that $\Gamma$ is piecewise analytic without cusps.
Then, for any closed set $E\subset\Omega$, there exists $n_0\in\mathbb{N}$, such that for $n\ge n_0$,
$p_n(z)$ \textit{has no zeros} on $E$.
\end{theorem}
\begin{theorem}
Assume that $\Gamma$ is piecewise analytic without cusps. Then,
$$
\lim_{n \to\infty}|p_n(z)|^{1/n}=|\Phi(z)|,\quad
z\in{\Omega}\setminus\{\infty\}.
$$
\end{theorem}

\bibliographystyle{amsplain}

\end{document}